\documentclass[12pt]{article}
\usepackage{latexsym,amsfonts,amssymb}

\usepackage{amsmath}
\usepackage{amsthm}
\usepackage{amssymb}
\usepackage{delarray}
\usepackage[dvips]{graphics}
\usepackage{epsfig}
\usepackage{color}
\usepackage{setspace}

\newtheorem{thm}{Theorem}[section]
\newtheorem{lemma}[thm]{Lemma}

\theoremstyle{definition}

\theoremstyle{remark}

\newcommand{\zm}{\noindent {\bf Proof: }}

\def\cd{\cdots}

\def\dim{\text{\rm dim}\,  }
\def\ch{\text{\rm ch}\,  }

\def\hom{\text{\rm Hom}\,  }

\def\al{\alpha}

\def\la{\lambda}

\def\g{\mathfrak g}

\begin{document}
\title{Rigidity of formal characters of  \\
        Lie algebras of type $A$
    }
\date{}
\author{Zhongguo Zhou  
}
\date{\small\it College of Science, Hohai University\\
    Nanjing,  210098, P.R.China.\\ {
{\rm e-mail: zhgzhou@hhu.edu.cn}\\}
}
\maketitle
\begin{abstract}
For a complex simple Lie algebra of type $A$, given a family of elements
$f_\la\in \mathbb Z[\Lambda],\la\in \Lambda^+,$
we show that $f_\la$ is just the formal character of the Weyl module $V(\la) $
if $f_\la$ satisfy two natural conditions.\\


\medskip
\noindent {\bf 2000 MR Subject Classification:} 17B10 \quad 20G05
\end{abstract}
\section{Introduction}

Let $\g_l$ be a complex simple Lie algebra of type $A_l$,
and let $V(\la)$ be the Weyl module. The formal character
$\ch_\lambda$ of $V(\la)$ is determined by the Weyl character
formula or by other methods, such as the Freudenthal formula.
The Littlewood-Richardson coefficient, $c_{\mu,\,\nu}^\la$
defined the multiplicity of $V(\la)$ in $V(\mu)\bigotimes V(\nu)$,
can be determined according to the formal character $\ch_\lambda.$
In this paper we shall consider the following problem:

\medskip
\textit{ Given a family of elements  $f_\la\in \mathbb
Z[\Lambda],\la\in \Lambda^+,$ when are they equal to $\ch_\la?$ }
\medskip

We prove that two natural conditions are enough to do so.
The one says it should contain right information about the
Weyl module for proper subalgebras of $\g_l,$  that is to say,
for all $\beta =\lambda-\mu, $ $m_\lambda(\mu)=n_\lambda(\mu), $
if $|\mbox{Supp}(\beta)|<l.$ The other says the number $n_{\mu,\,\nu}^\la $,
defined similarly to $c_{\mu,\,\nu}^\la$, should satisfy
$n_{\mu,\,\nu}^\lambda=n_{\lambda,\,-w_0\nu}^\mu.$

The motivation for this problem comes from the
modular representation theory of linear algebraic groups in
positive characteristic, especially the Lusztig conjecture.
In Lusztig's theory, the $\ch L(\lambda)$ is given by the
Kanzhdan-Lusztig polynomials and $\ch V(\mu)'s.$ The first condition
on $\ch L(\lambda)$ can be proved by an induction on rank. So we
want to find out some natural condition to ensure the conjecture.
Hence we consider the problem in complex field at first, the most simple case.

\section{notations}

Let $\g_l$ be a complex simple Lie algebra of type $A_l$,
and let $\Delta=\{\al_1, \al_2,\cd, \al_l\}$ be the set of simple roots.
Let $\Lambda$ be the set of weights, and $\omega_1, \omega_2,
\cdots, \omega_l$ the set of fundamental dominant weights. Then the
set of dominant weights is denoted by $\Lambda^+.$ Let $\ch_\lambda$,
$\lambda\in \Lambda^+$, be the formal character of the Weyl module
$V(\lambda)$, they form a free $\mathbb Z$-module of the commutative
ring $\mathbb Z[\Lambda]$, with base $\{e(\lambda),\lambda\in\Lambda\},$
and multiplication $e(\lambda)*e(\mu)=e(\lambda+\mu).$  Let $W$ be the
Weyl group, action on  $\mathbb Z[\Lambda]$ naturally as $\sigma e(\lambda)
=e(\sigma\lambda).$ Set
$$
\mathbb Z[\Lambda]^W=\{f\in \mathbb Z[\Lambda]\,|\, wf=f, \;\;w\in W\}.
$$
Let $W_\lambda $ be the $W$-orbit of $\lambda$ and $h(\lambda )=\sum_{x\in
W_\lambda} e(x).$ Let $\Pi(\lambda)$ be the set of saturated weights of
weight $\lambda$ and $\Pi^+(\lambda)=\Pi(\lambda)\cap \Lambda^+.$
It is well known that
$$ \ch_\lambda=\sum_{x\in\Pi^(\lambda)}m_\lambda(x)e(x)=\sum_{\mu\in
\Pi^+(\lambda)}m_\lambda(\mu)h(\mu), \;\;\lambda\in \Lambda^+,$$
forms a basis of $\mathbb Z[\Lambda]^W.$

Recall that $w_0$ is the longest element  of $W$ with the action
$w_0\omega_i=-\omega_{l-i}.$  Let $c_{\mu,\,\nu}^\la$  be the
Littlewood-Richardson coefficient.  According to the complete
reducibility of $\g_l$-modules and
$$
\dim\hom(U\otimes V, W)=\dim\hom(U ,W\otimes V^*),
$$
we have
$$
[V(\mu)\otimes V(\la): V(\nu)]=[V(\nu)\otimes V(\la)^*: V(\mu)]=[V(\nu)\otimes V(-w_0\la): V(\mu)].
$$
Hence
$$
c_{\mu,\,\la}^\nu=c_{\nu,\,-w_0\la}^\mu.
$$

\section{main results}

{\bf 3.1.} Let $\beta=\sum_{i=1}^{l}k_i\alpha_i,$ define
$$\mbox{Supp}(\beta)=\{\alpha_i\,| \,k_i\neq 0\}.$$
For $\lambda \in \Lambda^+,$ set $f_{\lambda}
=\sum_{\mu\in \Pi^+(\lambda)}n_\lambda(\mu)h(\mu)
=\sum_{x\in \Pi(\lambda)}n_\lambda(x)e(x)
\in \mathbb Z[\Lambda]^W$  satisfied $n_\lambda(\lambda)=1, n_\lambda(x)=0,$ if
$x\notin \Pi(\lambda).$ Then $f_{\lambda}$ is also a basis of $\mathbb Z[\Lambda]^W.$
Hence there exists unique $n_{\mu,\,\nu,}^\lambda $ such that
$$
f_{\mu}* f_{\nu}=\sum_{\lambda\in\Pi^+(\mu+\nu)}n_{\mu,\,\nu}^\lambda f_{\lambda}.
$$
By the definition of $f_{\lambda}'s,$ we have
$$n_{\mu,\,\nu}^{\mu+\nu}=1; n_{\mu,\,\nu}^\lambda=n_{\nu,\,\mu}^\lambda.$$

\begin{thm}\label{th1}
If these $f_{\lambda}$ satisfy the following two conditions:

(1) for all $\beta =\lambda-\mu, m_\lambda(\mu)=n_\lambda(\mu)$, if $|{\rm{Supp}}(\beta)|<l.$

(2) $n_{\mu,\,\nu}^\lambda=n_{\lambda,\,-w_0\nu}^\mu.$

\noindent Then $f_{\lambda}=ch_\lambda, n_{\mu,\,\nu}^\lambda=c_{\mu,\,\nu}^\lambda.$
\end{thm}
The two conditions in theorem \ref{th1} are very natural. The first
one comes from restricting $V(\lambda)$ to subalgebra with smaller
rank and the second one is because the tensor-hom adjunction. So it
is surprised that $ch_\lambda, n_{\mu,\,\nu}^\lambda$ can be determined
completely only by two conditions. Hence we call this property as the
rigidity of formal characters or representation ring of Lie algebra.

Firstly, we need a lemma.

\begin{lemma}
For $f_{\lambda},$ the numbers $n_{\mu,\,\nu}^\lambda$ are uniquely
determined  by $n_\lambda(\mu),$ and vice versa.
\end{lemma}

\zm Because $f_{\mu}* f_{\nu}\in \mathbb Z[\Lambda]^W$ and $f_\lambda,
\lambda\in \Lambda^+$, form the basis of $\mathbb Z[\Lambda]^W,$ so
 $f_{\mu}* f_{\nu}$ can be written in terms of a linear
combination of base elements $f_{\lambda}, \lambda\in \Lambda^+$, i.e.
$$f_{\mu}* f_{\nu}=\sum_{\lambda\in\Pi^+(\mu+\nu)}n_{\mu,\,\nu}^\lambda f_{\lambda}.$$
Therefore all these $n_{\mu,\,\nu}^\lambda$ are uniquely determined.

Now suppose that we know $n_{\mu,\,\nu}^\lambda$ already, and will determine
$n_\lambda(\mu)$ recursively according to a suitable partial order in $\Lambda^+$
as follow.

Let $\lambda,\mu\in \Lambda^+,$ define $\mu< \lambda$ if $\mu\prec \lambda $  or
$ \lambda-\mu\in   \Lambda^+.$

For the fundamental weights $\omega_i$ and $0,$ their saturated weight set $\Pi(\omega_i)$
only contains one dominate weight. So
$$f_{\omega_i}=h(\omega_i), \;\; f_{0}=h(0)=e(0).$$
Suppose that $\lambda\in \Lambda^+$ does not be a fundamental
dominate weight or 0, then there exist $\mu,\nu\in\Lambda^+,$ such that
$$\lambda=\mu+\nu, \; \mu< \lambda,\; \nu< \lambda.$$
Then
$$f_{\mu}* f_{\nu}=f_{\lambda}+\sum_{s\in \Pi^+(\lambda),\, s\neq \lambda}n_{\mu,\,\nu}^s f_{s},$$
i.e.
$$f_{\lambda}=f_{\mu}* f_{\nu}-\sum_{x\in \Pi^+(\lambda),\, s\neq \lambda}n_{\mu,\,\nu}^s f_{s}.$$
Therefore we can get a precise expression for $f_{\lambda}$ because all
$f_{\mu}, f_{\nu}, f_{s}, n_{\mu,\,\nu}^s$ can be determined recusively
as $\mu< \lambda, \nu< \lambda ,s< \lambda$. The lemma is proved.

{\bf Remark}: This lemma doesn't work for type $B_l.$  For example, for $B_2,$
we can't determine $(1,0)_{(0,0)}$ by using the above method because there is no
way to decompose $(1,0)$ into a sum of two non-trivial dominate weights.

{\bf3.2.} From the equation $f_{\lambda}=f_{\mu}* f_{\nu}-\sum_{x\in \Pi^+(\lambda),
\, s\neq \lambda}n_{\mu,\,\nu}^s f_{s}$, we have more precise expression.
\begin{eqnarray*}
f_{\lambda}\!\!\!\!\!\!&&=\sum_{t\in \Pi(\lambda)}n_\lambda(t)e(t)\\
&&=
\sum_{y\in \Pi(\mu)}n_\mu(y)e(y)*
\sum_{z\in \Pi(\nu)}n_\nu(z)e(z)
-\sum_{s\in \Pi^+(\lambda),s\neq \lambda}n_{\mu,\,\nu}^s
\sum_{x\in \Pi(s)}n_s(x)e(x)\\
&&=
\sum_{y\in \Pi(\mu),z\in \Pi(\nu)}n_\mu(y)n_\nu(z)e(y+z)
-\sum_{s,x\in \Pi(\lambda),s\neq \lambda,}n_{\mu,\,\nu}^s n_s(x)e(x)\\
&&=
\sum_{t\in\Pi(\lambda)}\left(
\sum_{y\in \Pi(\mu),z\in \Pi(\nu),y+z=t}n_\mu(y)n_\nu(z)\right)e(t)-\\
&&-\sum_{t\in \Pi(\lambda)}\left(\sum_{s\in \Pi(\lambda),s\neq \lambda,}n_{\mu,\,\nu}^s n_s(t)\right)e(t).
\end{eqnarray*}
Hence
\begin{eqnarray}{\label{gsi}}
n_\lambda(t)=
\sum_{y\in \Pi(\mu),z\in \Pi(\nu),y+z=t}n_\mu(y)n_\nu(z)
-\sum_{s\in \Pi(\lambda),s\neq \lambda}n_{\mu,\,\nu}^s n_s(t),
\end{eqnarray}
and then
\begin{eqnarray}
n_{\mu,\,\nu}^t
=
\sum_{y\in \Pi(\mu),z\in \Pi(\nu),y+z=t}n_\mu(y)n_\nu(z)
-\sum_{s\in \Pi(\lambda),s\neq t,}n_{\mu,\,\nu}^s n_s(t).
\end{eqnarray}

Note from the two formulas that $n_\lambda(t)$ is only depended on those
numbers
\begin{eqnarray*}
& n_\mu(y), \quad &  \mu-y\preceq \lambda-t, \mu<\lambda;\\
 & n_\nu(z) ,\quad &   \nu-z\preceq\lambda-t,\nu< \lambda; \\
 & n_s(t),    n_{\mu,\,\nu}^s,  \quad   &  t\preceq s\prec \lambda
\end{eqnarray*}
and  $n_{\mu,\,\nu}^t$ is only depended on those
numbers
\begin{eqnarray*}
& n_\mu(y), \quad &  \mu-y\preceq \lambda-t, \mu<\lambda;\\
 & n_\nu(z) ,\quad &   \nu-z\preceq\lambda-t,\nu< \lambda; \\
 & n_s(t),    n_{\mu,\,\nu}^s,  \quad   &  t\prec s\preceq \lambda.
\end{eqnarray*}
We can substitute $n_{\mu,\,\nu}^s,s\succ t$, in (2) for those in (1),
and obtain
\begin{eqnarray}\label{eqn3}
n_\lambda(t)=n_{\mu,\,\nu}^t +g(n_\mu(y),n_\nu(z),n_s(x)),
\end{eqnarray}
where $g(n_\mu(y),n_\nu(z),n_s(x)) $ is determined by
$ n_\mu(y), n_\nu(z), n_s(x),t\preceq x,s\preceq \lambda$ and $ s-x\prec\lambda-t.$

\medskip
{\bf 3.3.}  We now prove the theorem by induction on $\Lambda^+$ with
the partial order $\lq\lq<"$ and on $\Pi(\lambda)$ with the partial order
$\lq\lq\prec".$

It is only need to prove $n_\lambda(\mu)=m_\lambda(\mu), \lambda\in \Lambda^+,\mu\in\Pi^+(\lambda).$

Firstly,  if $\lambda=\omega_i$ or  $0,$ their saturated weight set  $\Pi(\lambda)$
only contains one dominate weight. So the theorem holds by the definition of $f_\lambda.$

Suppose that $\lambda\in \Lambda^+, $ not be fundamental weights $\omega_i$ or 0. Then $n_\lambda(\lambda)=1=m_\lambda(\lambda).$ Let $\mu\in \Pi^+(\lambda),\beta=\lambda-\mu.$
Consider the three cases as follows:

(1)\,  when $|\rm{Supp}(\beta)|<l,$ we have
$n_\lambda(\mu)=m_\lambda(\mu)$ by the first condition in theorem.

(2)\,  when $|\rm{Supp}(\beta)|=l $ and $\beta=\alpha_1+\alpha_2+\cdots+\alpha_l=\omega_1+\omega_l,$
we have $\mu=\lambda-\beta=\lambda-\omega_1-\omega_l,$ thus $\lambda_1\geq 1,\lambda_l\geq 1.$
Moreover, $\lambda=(\lambda-\omega_1)+\omega_1, \lambda-\omega_1\in \Pi^+(\lambda)$.

By the second condition of theorem, we have
$$n_{\omega_1,\,\lambda-\omega_1}^\mu=n_{\omega_1,\,\lambda-\omega_1}^{\lambda-\omega_1-\omega_l}
=n_{\lambda-\omega_1-\omega_l,\,-w_0\omega_1}^{\lambda-\omega_1}
=n_{\lambda-\omega_1-\omega_l,\,\omega_l}^{\lambda-\omega_1}=1.$$
The last equation holds  because
$$
{\lambda-\omega_1-\omega_l+\omega_l}-({\lambda-\omega_1})=0.
$$
We also have
$$
    c_{\omega_1,\,\lambda-\omega_1}^\mu
   =c_{\omega_1,\,\lambda-\omega_1}^{\lambda-\omega_1-\omega_l}
   =c_{\lambda-\omega_1-\omega_l,\,-w_0\omega_1}^{\lambda-\omega_1}
   =c_{\lambda-\omega_1-\omega_l,\,\omega_l}^{\lambda-\omega_1}
   =1.
$$
Therefore,
$$
n_{\omega_1,\,\lambda-\omega_1}^\mu
=c_{\omega_1,\,\lambda-\omega_1}^\mu.
$$
For $\omega_1<\lambda, \,\lambda-\omega_1<\lambda,$ according to equation (\ref{eqn3}) and induction
hypothesis,  we have
\begin{eqnarray*}
&&n_\lambda(\mu)=n_{\omega_1,\,\lambda-\omega_1}^\mu +g(n_{\omega_1}(y),n_{\lambda-\omega_1}(z),n_s(x))\\
&&=c_{\omega_1,\,\lambda-\omega_1}^\mu +g(c_{\omega_1}(y),c_{\lambda-\omega_1}(z),c_s(x))=c_\lambda(\mu).
\end{eqnarray*}
(3)\,  When $|\rm{Supp}(\beta)|=l $ and $\beta=\alpha_1+\alpha_2+\cdots+\alpha_l+\beta_1,\beta_1\succ 0,$
we have $\mu=\lambda-\omega_1-\omega_l-\beta_1.$ If  $\lambda_i> 0, \lambda_j=0 $ and $j<i,\,$ let $\nu=\lambda-\omega_i, $ then $\lambda=\nu+\omega_i,$ and
$$
n_{\omega_i,\,\nu}^\mu
=n_{\omega_i,\,\lambda-\omega_i}^{\lambda-\omega_1-\omega_l-\beta_1}
=n_{\lambda-\omega_1-\omega_l-\beta_1,\, -w_0\omega_i}^{\lambda-\omega_i}
=n_{\lambda-\omega_1-\omega_l-\beta_1,\, \omega_{l-i}}^{\lambda-\omega_i}
=n_{\mu,\, \omega_{l-i}}^{\nu}.
$$
For
\begin{eqnarray*}
&&\mu+\omega_{l-i}-{\nu}\\
&&={\lambda-\omega_1-\omega_l-\beta_1+ \omega_{l-i}}-{\lambda+\omega_i}\\
&&=-(\omega_1+\omega_l)-(\omega_{l-i}+\omega_i)-\beta_1\\
&&=-\sum_{k=1}^l\alpha_k+\sum_{k=1}^{min(l-i,\,i)}\sum_{j=k}^{l-k}k\alpha_k-\beta_1,
\end{eqnarray*}
let $\mu+\omega_{l-i}-{\nu}=\sum_{j=1}^lk_j\alpha_j,$ then $k_1\leq 0.$

(i) \, if $k_1< 0,$  then $\nu\notin \Pi(\mu+\omega_{l-i}),$ we have  $n_{\mu,\, \omega_{l-i}}^{\nu}=0
=c_{\mu,\, \omega_{l-i}}^{\nu}.$

(ii)\, if $k_1= 0,$  then $|\rm{Supp}(\mu+\omega_{l-i}-{\nu})|<l,$ we have by the first condition
$$
n_{\mu,\, \omega_{l-i}}^{\nu}
=c_{\mu,\, \omega_{l-i}}^{\nu}.
$$
Moreover,
$$
c_{\omega_i,\,\nu}^\mu
=c_{\omega_i,\,\lambda-\omega_i}^{\lambda-\omega_1-\omega_l-\beta_1}
=c_{\lambda-\omega_1-\omega_l-\beta_1,\, -w_0\omega_i}^{\lambda-\omega_i}
=c_{\lambda-\omega_1-\omega_l-\beta_1,\, \omega_{l-i}}^{\lambda-\omega_i}
=c_{\mu,\, \omega_{l-i}}^{\nu}.
$$
Therefore, we have also
$$
n_{\omega_i,\,\nu}^\mu
=c_{\omega_i,\,\nu}^\mu.
$$
Similarly, by the induction hypothesis and equation (\ref{eqn3}) we have
$$n_\lambda(\mu)=c_\lambda(\mu).$$
The theorem is proved.

\medskip

{\bf 3.4.} One of the conditions in theorem is about $m_\lambda(\mu),$ and the other is
about $n_{\mu,\,\nu}^\lambda.$ There exists difficulty to generalize the theorem to Lie
algebras of other type. Only the two conditions in theorem are not enough to determined
$m_\lambda(\mu).$ We need to find out a suitable condition. As said before, the motivation
for this problem in this paper comes from the modular representation theory of linear
algebraic groups in positive characteristic. It is more difficulty to find out natural
conditions in this case, for the second condition does not hold in general. We will consider
these problems in the future.

\medskip

\noindent{\bf Acknowledgement:} This work was
supported  by the Fundamental Research Funds for the Central
Universities 2009B26914 and 2010B09714.

\makeatletter
\def\@biblabel#1{#1.}
\makeatother

\end{document}